\documentclass[11pt]{article}
\usepackage{graphicx}
\usepackage{amsmath}
\usepackage{amssymb}
\usepackage{amsthm}
\usepackage{epsfig}
\usepackage{latexsym} 
\usepackage{color}

\newtheorem{theorem}{Theorem}[section]
\newtheorem{lemma}[theorem]{Lemma}
\newtheorem{corollary}[theorem]{Corollary}
\newtheorem{conjecture}[theorem]{Conjecture}

\newtheorem{proposition}[theorem]{Proposition}

\def\PSL{\mbox{\rm{PSL}}} 

\def\Isom{\mbox{\rm{Isom}}} 
\def\SL{\mbox{\rm{SL}}} 
\def\PGL{\mbox{\rm{PGL}}} 
\def\tr{\mbox{\rm{tr}}} 
\def\orb{\mbox{\rm{orb}}} 
\def\V{V}
\title{Commensurability classes of 2-bridge knot complements} 
\author{A. W. Reid\thanks{This work was 
partially supported by the N. S. F.} ~ \&~ G. S. Walsh}

\date{}
\begin{document}
\maketitle
\abstract{We show that a hyperbolic 2-bridge knot complement is the unique knot complement in its commensurability class. We also discuss constructions of
commensurable hyperbolic knot complements and put forth a conjecture 
on the number of hyperbolic knot complements in a commensurability class.}

\section{Introduction}
Recall that two hyperbolic 3-manifolds $M_1={\bf H}^3/\Gamma_1$ and
$M_2={\bf H}^3/\Gamma_2$ are {\it commensurable} if they have 
homeomorphic finite sheeted covering spaces.  In terms of the groups,
this is equivalent to $\Gamma_1$ and some conjugate of $\Gamma_2$ in
$\Isom({\bf H}^3)$ having a common finite index subgroup.  Proving
that two hyperbolic 3-manifolds are commensurable (or not
commensurable) is in general a difficult problem.  The most useful
techniques at present are algebraic; for example the invariant
trace-field (see \cite{MR} Chapter 3).

In this paper we investigate commensurability of hyperbolic knot
complements in $S^3$.  Our main result is:
 
\begin{theorem} 
\label{main}
Let $K$ be a hyperbolic 2-bridge knot.  
Then $S^3\setminus K$ is not commensurable with the complement of any other 
knot in $S^3$. 
\end{theorem}  
  
Previous work in this direction was done in \cite{Re1}, where it is
shown that the figure-eight knot is the only knot in $S^3$ whose
complement is arithmetic.  In addition, it is known that no two
hyperbolic twist knot complements are commensurable \cite{HS}.  On
the other hand, there are hyperbolic knot complements in $S^3$ with
more than one knot complement in the commensurability class. A detailed
discussion of this is given in \S 5.

A corollary of our main theorem is the following result
which is a direct consequence of Corollary 1.4 of \cite{Sc}.

\begin{corollary}
\label{QI}
Let $K$ be a hyperbolic 2-bridge knot in $S^3$ and $K'$ any knot in
$S^3$.  If $\pi_1(S^3\setminus
K)$ and $\pi_1(S^3\setminus K')$ are quasi-isometric, then $K$ and
$K'$ are equivalent.
\end{corollary}
 
\medskip

\noindent{\bf Acknowledgements:}~The authors wish to thank Michel Bolieau,  Cameron Gordon and Walter Neumann for some useful
conversations on matters related to this work. The first author also
wishes to thank the Institute for Advanced Study, and the Universidade
Federal do Cear\'a, Fortaleza for their hospitality during the
preparation of this work.

\section{Preliminaries} \label{background}

We begin by recalling some terminology and results that will be
needed. Henceforth, unless otherwise stated {\em knot complement} will
always refer to a knot complement in $S^3$.

\subsection{Hidden Symmetries} 

Let $\Gamma$ be a Kleinian group of finite co-volume. The {\it commensurator}
of $\Gamma$ is the group

$$C(\Gamma) = \lbrace g \in \Isom({\bf H}^3): [\Gamma: \Gamma \cap g^{-1} \Gamma g] <\infty\}$$

We denote by $C^+(\Gamma)$ the subgroup of $C(\Gamma)$ of index at most 2 that consists of orientation-preserving isometries.  It is a fundamental result of Margulis \cite{Ma}
that $C^+(\Gamma)$
is a Kleinian group if and only if $\Gamma$ is non-arithmetic, and moreover, in this case,  $C^+(\Gamma)$ is the unique maximal element in the $\PSL(2,\mathbb{ C})$ commensurability class of $\Gamma$. 

Note that the normalizer of $\Gamma$ in $\PSL(2,\mathbb{ C})$, which
we shall denote by $N(\Gamma)$, is a subgroup of $C^+(\Gamma)$.
In the case when $\Gamma$ corresponds to the
faithful discrete representation of $\pi_1(S^3\setminus K)$ ($K$
distinct from the figure-eight knot) it is often the case that
$N(\Gamma) = C^+(\Gamma)$.  Before we give a more precise
discussion of this, we recall some of \cite{NR}. Henceforth, any knot
will be assumed hyperbolic and distinct from the figure-eight knot.

Assume that
$S^3\setminus K = {\bf H}^3/\Gamma_K$. $K$ is said to have {\em hidden
symmetries} if $C^+(\Gamma_K)$ properly contains $N(K)=N(\Gamma_K)$.  
We will make use of the following result 
 from \cite{NR}, which requires one more piece of terminology.

Let $S^2(2,4,4)$, $S^2(2,3,6)$ and $S^2(3,3,3)$ denote the
2-dimensional orbifolds which are 2-spheres with 3 cone points and
cone angles $(\pi, \pi/2,\pi/2)$, $(\pi, 2\pi/3, \pi/3)$ and $(2\pi/3,
2\pi/3, 2\pi/3)$ respectively. These are called {\em Euclidean
  turnovers}.  In addition, we let $S^2(2,2,2,2)$ denote the 
2-dimensional orbifold which is 2-sphere with 4 cone points all of 
cone angle $\pi$.

If $X$ is an orientable, non-compact finite volume
hyperbolic 3-orbifold, then a cusp of $X$ has the form $Q \times [0,
\infty)$,  where $Q$ is a Euclidean orbifold. The cusp is said to be {\em rigid} if $Q$
is a Euclidean turnover.

\begin{proposition}  \cite[Proposition 9.1]{NR}
\label{nr}
The following are equivalent for a hyperbolic knot $K$
other than the figure eight knot complement:\\[\baselineskip]
\noindent (i)~it has hidden symmetries;\\[\baselineskip]
\noindent (ii)~the orientable commensurator quotient ${\bf H}^3/C^+(\Gamma_K)$ has 
a rigid cusp;\\[\baselineskip]
\noindent (iii)~$S^3\setminus K$ non-normally covers some orbifold.
\end{proposition}
 
\medskip

\noindent{\bf Notation:}~In what follows, we shall let
$Q_K={\bf H}^3/C^+(\Gamma_K)$.

\subsection{Cusp Fields and Trace Fields}

Recall that if $\Gamma$ is a Kleinian group of finite co-volume, the
trace field is a finite extension of $\mathbb{ Q}$.  Furthermore,   the
invariant trace-field of $\Gamma$, $k\Gamma =
\mathbb{Q}(\tr(\gamma^2) : \gamma \in \Gamma)$, is a subfield of the
trace-field that is an invariant of the commensurability class 
\cite[Chapter 3]{MR}.  When $\Gamma_K$ is the faithful discrete
representation of $\pi_1(S^3\setminus K)$, it is shown in
\cite{Re} that the invariant trace-field coincides with the
trace-field.  This also holds when the Kleinian group is generated by
parabolic elements \cite[Lemma 1]{Re1}.

 For
convenience we will often abuse the distinction between $\PSL$ and
$\SL$ and simply work with matrices. If ${\bf H}^3/\Gamma$ is a 1-cusped hyperbolic 3-manifold, we can conjugate the
peripheral subgroup to be $<\begin{pmatrix} 1 & 1 \\ 0 & 1
\end{pmatrix}, \begin{pmatrix} 1 & g \\ 0 & 1 \end{pmatrix}>$. It is easily shown that $g\in
k\Gamma$ (see for example \cite[Proposition 2.7]{NR}) and $g$ is
referred to as the {\em cusp parameter} of $\Gamma$.   In
  the natural identification of the Teichm\"uller space of the torus
  with the upper half-plane, $g$ is the shape of the torus. The field
$\mathbb{Q}(g)$ is called the {\em cusp field}, which is a subfield of
$k\Gamma$.
                                                           
Of relevance to us is that there are constraints on cusp parameters of tori that cover turnovers.  The rigid orbifolds 
$S^2(2,4,4), S^2(2,3,6)$ and $S^2(3,3,3)$
have orbifold groups that are extensions of $\mathbb{Z} \oplus \mathbb{Z}$
by elements of orders $4$, $6$ and $3$ respectively. The maximal $\mathbb{Z} \oplus \mathbb{Z}$
subgroup in these cases can be conjugated to be:
$$<\begin{pmatrix} 1 & 1 \\ 0 & 1 \end{pmatrix},
                    \begin{pmatrix} 1 & \epsilon \\ 0 & 1\end{pmatrix}>,$$
where $\epsilon = i$ in the case of $S^2(2,4,4)$ and 
$\epsilon =(-1+\sqrt{-3})/2$ (which we shall denote by $\omega$ in what follows) otherwise.  
This discussion together with Proposition \ref{nr}, yields the following corollary. 

\begin{corollary}
\label{cuspfield}
Let $K$ be a hyperbolic knot which has
hidden symmetries. Then the cusp parameter of $S^3\setminus K $ lies 
in $\mathbb{Q}(i)$ (when the turnover is $S^2(2,4,4)$) or in $\mathbb{Q}(\sqrt{-3})$ (when the
turnover is $S^2(2,3,6)$ or $S^2(3,3,3)$).
 \end{corollary}

%Since $H_1(M \setminus K) = \mathbb{Z}$, any such covering must be cyclic with covering group generated by $\tau$.   Fill in $M' \setminus  K'$ to obtain the homotopy 3-sphere $M'$.  Then since $K'$ is non-trivial, the induced action of $\tau$ on $M'$ is free, by the solution to the Smith conjecture, \cite{MB}. Therefore $M' / <\tau>$ is a  manifold with cyclic fundamental group obtained by surgery on $M\setminus K$. Now suppose there is a cyclic filling of $M \setminus K$ along a non-meridional slope $\alpha$.  Denote this by $M(\alpha)$  Let $M'$ be the  universal cover of $M(\alpha)$.  Denote the pre-image of $K$ by $K'$.  Then $M' \setminus K'$ is a cyclic cover of $M \setminus K$, and $H_1(M \setminus K, \mathbb{Z}) = \mathbb{Z}$. Thus $K'$ must have one  component.  Note that $\alpha$ lifts to $M' \setminus K'$. \qed

\subsection{2-bridge knots}

It will be convenient to recall some facts about 2-bridge knots that
we shall make use of.  In particular the work of Riley (\cite{Ri1},
\cite{Ri2}) is heavily used. Thus throughout this section let $K$ be a
hyperbolic 2-bridge knot.

A 2-bridge knot $K$ has a normal form $(p,q)$ where 
$p$ and $q$ are odd integers and are determined by the lens space $L(p,q)$
that is the double cover of $S^3$ branched over $K$. 
The fundamental group of a 2-bridge knot complement has a presentation
of the form $\pi_1(S^3 \setminus K) = <x_1,x_2:r>$ where $x_1$ and $x_2$ are
meridians and the relation $r$ has the form $wx_1w^{-1} = x_2$ for some
word $w$ in $x_1$ and $x_2$. The exponents of the $x_i$ in the word
$w$ are all $\pm 1$, and are
determined by the 2-bridge normal form of $K$.

Let $\bf F$ be a field and fix an algebraic closure ${\overline{\bf
    F}}$. A homomorphism $\rho:\pi_1(S^3\setminus K)\rightarrow
\PSL(2,{\bf F})$ is called a {\em parabolic representation} (or simply
{\em p-rep} for short) if $\rho(x_1)$ (and hence $\rho(x_2)$) is a
parabolic element; i.e. conjugate in $\PSL(2,\overline{\bf F})$ to the
element $\begin{pmatrix} 1 & 1 \\ 0 & 1 \end{pmatrix}$.

If we conjugate so as to consider a p-rep normalized so that
$$\rho(x_1) = \begin{pmatrix} 1 & 1 \\ 0 & 1 \end{pmatrix}~\hbox{and}~\rho(x_2) = \begin{pmatrix} 1 & 0 \\ y & 1 \end{pmatrix},$$
then Riley shows that $y$ satisfies a certain polynomial $\Lambda_K(y)$ with 
leading coefficient and constant term equal to 1 (\cite{Ri1} Theorem 2).  
We shall say that the above p-rep is in {\em standard form}.

In the case that ${\bf F}=\mathbb{ C}$, $\Lambda_K$ is a polynomial
with integral coefficients and some root of the p-rep polynomial
corresponds to the faithful discrete representation of $\pi_1(S^3
\setminus K)$ into $\PSL(2,\mathbb{ C})$. In this case, $\Lambda_K(y)$
is called the {\em p-rep polynomial of $K$}.

In \cite{Ri1}, Riley also describes the image of a longitude $\ell$
for $x_1$ for a p-rep in standard form; namely a matrix
$\begin{pmatrix} 1 & g \\ 0 & 1 \end{pmatrix}$, where $g = 2 g_0$ and
$g_0$ is an algebraic integer in the field $\mathbb{ Q}(y)$.

Also important in what follows is Riley's work on p-reps when $\bf F$
has characteristic $2$. Here the images of the
meridians are elements of order $2$, and so the image of a p-rep is
therefore a dihedral group.  This dihedral group is necessarily finite
since a knot group cannot surject the infinite dihedral group. 
In addition, since the image groups are non-cyclic,
the dihedral groups considered are never of order $2$.
Hence we exclude this case from further comment.
Riley proves the following result which will be useful for us (see
\cite{Ri1} Theorem 3).  

\begin{theorem}
\label{norepeated}
The p-rep polynomial $\Lambda_K(y)$ has no repeated factors modulo $2$, and so no repeated factors
in $\mathbb{ Z}[y]$.
\end{theorem} 

This result allows us to prove Proposition \ref{no2} below (which Riley noticed in
\cite{Ri2}). For completeness we give a proof. We first record 
the following standard
facts about polynomials reduced modulo primes (for example see \cite[Proposition 3.8.1 and Theorem 3.8.2]{Ko}).

\begin{theorem}
\label{numbertheory}
Let $f(x)\in\mathbb{Z}[x]$ be an irreducible monic polynomial,
$\alpha$ a root and $k=\mathbb{Q}(\alpha)$ with ring of integers
$R_k$.  Let $d_k$ denote the discrminiant of $k$, and $\Delta(\alpha)$
the discriminant of $f$. Let $p$ be a rational prime and
$\bar{f}$ the reduction of $f$ modulo $p$. Then,

\smallskip

\item{(i)}~$\bar{f}$ decomposes into distinct irreducible factors if and
only if $p$ does not divide $\Delta(\alpha)$.

\smallskip

\item{(ii)}~Suppose that $p$ does not divide $\Delta(\alpha)d_k^{-1}$ and 
$\bar{f} = \bar{f}_1^{e_1} \ldots \bar{f}_g^{e_g}.$
Then $pR_k = {\cal P}_1^{e_1} \ldots {\cal P}_g^{e_g}$
is the factorization into prime powers.
\end{theorem}

 \begin{proposition} \label{no2} Let $K$ be a hyperbolic 2-bridge knot with trace-field $k$.  Then 
   $\mathbb{ Q}(i)$ is not a subfield of $k$.
 \end{proposition}

\noindent{\bf Proof:}~We shall show that $2$ does not divide the discriminant of $k$. Since
the discriminant of $\mathbb{ Q}(i)$ is $-4$, standard facts about
the behavior of discriminant in extensions of number fields precludes
$\mathbb{ Q}(i)$ from being a subfield of $k$ (see \cite{Ko} for
example).

Let $\rho$ denote the p-rep corresponding to the faithful
discrete representation conjugated to be in standard form, and
$\Lambda_0(y)$ be the irreducible factor of $\Lambda_K(y)$ which gives
the representation corresponding to the complete structure. We denote
the image group by $\Gamma_K$. Therefore $k=k\Gamma_K = \mathbb{Q}(y)$
for some root $y$ of $\Lambda_0(y)$.  By Theorem \ref{norepeated}
$\Lambda_K(y)$ has distinct factors modulo 2, and so $\Lambda_0(y)$
has distinct factors modulo 2. Thus, Theorem \ref{numbertheory}(i)
shows that $2$ does not divide the discriminant $\Delta(y)$ of
$\Lambda_0(y)$.  Since the discriminant of $\mathbb{Q}(y)$ divides
$\Delta(y)$ (see \cite{Ko} Chapter 3)
it follows that 2 does 
not divide the discriminant of $\mathbb{Q}(y)$. \qed\\[\baselineskip]
Proposition \ref{no2} shows that the
cusp field of a hyperbolic 2-bridge knot is not $\mathbb{Q}(i)$. Hence
if $K$ has hidden symmetries, Corollary \ref{cuspfield} shows that the
cusp of the orbifold $Q_K$ is either $S^2(3,3,3)$ or $S^2(2,3,6)$.
In addition, notice that the element $\mu=\begin{pmatrix} i & 0 \\ 0 & -i
\end{pmatrix}$ normalizes the p-rep  $\rho$ in standard form.
%$$\rho(x_1) = \begin{pmatrix} 1 & 1 \\ 0 & 1 \end{pmatrix}
%~\hbox{and}~\rho(y_1) = 
%\begin{pmatrix} 1 & 0 \\ y & 1 \end{pmatrix},$$$
Hence $\mu\in N(K)<C^+(\Gamma_K)$, and we deduce:
 
\begin{corollary} \label{noi} 
If $K$ is a hyperbolic 2-bridge knot with hidden symmetries the cusp of the orbifold $Q_K$ is  $S^2(2,3,6)$ and the cusp field is $\mathbb{ Q}(\sqrt{-3})$.\qed
\end{corollary} 

\subsection{Symmetry groups of 2-bridge knots}

The following discussion is well-known, but will be convenient
for us to include.

We first state a result (originally due to Conway) about
the symmetry groups of 2-bridge knot complements (a proof can be found in
\cite{Sa} for example). \\[\baselineskip]
\noindent{\bf Notation:}~
Throughout the paper, we let $\V$ denote the group $\mathbb{Z}/ 2
\mathbb{Z} \oplus \mathbb{Z}/ 2 \mathbb{Z}$ and $D_n$ the dihedral
group with $2n$ elements.

\begin{theorem}
\label{symms}
Let $K$ be a hyperbolic 2-bridge knot. Then 
$\Isom^+(S^3\setminus K) = N(K)/\Gamma_K$ is either
$\V$ or $D_4$. 
In both cases every non-trivial element acts non-freely (ie with non-empty fixed point set).
\end{theorem}

We need some additional information about the quotient orbifold 
${\bf H}^3/N(K)$ when $K$ is 2-bridge.

\begin{lemma}
\label{oneZ4}
 Let $K$ be a hyperbolic 2-bridge knot and $Q = {\bf H}^3/N(K)$. Then 

\smallskip
\noindent (i)~$N(K)$ is generated by elements of order $2$. 

\smallskip

\noindent (ii)~There is a unique two-fold cover $Q'={\bf H}^3/\Gamma$ 
of $Q$ with a torus cusp. 
Furthermore, all torsion elements of $\Gamma^{ab}$ have order 2. 
\end{lemma}

\noindent{\bf Proof:}~We begin with a discussion of 
a particular subgroup (isomorphic to V) of the symmetry group
of any hyperbolic 2-bridge knot complement.  It is well-known that 
$(S^3\setminus K)/V$ is an orbifold whose orbifold group is generated
by involutions. However, it will be convenient for what follows
to recall a geometric description of this.

\begin{figure}[h]
\begin{center}
\input{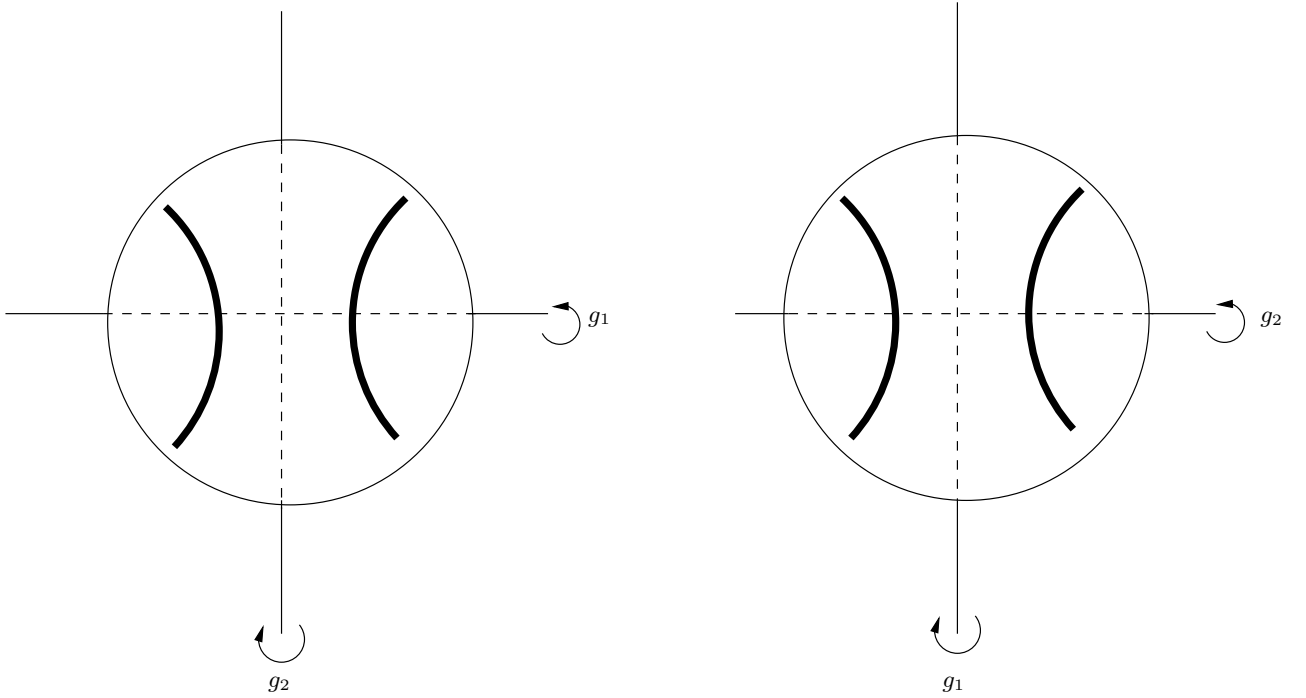}
\end{center}
\caption{The action of $V$ on a 2-bridge knot complement}
\label{Vsymmetries}
\end{figure}

The complement of a two-bridge knot in $S^3$ is the union of two
three-balls $B_1$ and $B_2$ with two unknotted arcs deleted, see
Figure \ref{Vsymmetries}.  We take the two order two symmetries $g_1$
and $g_2$ which pointwise fix the centers of the arcs in $B_1$ and
$B_2$ respectively to be the generators of $\V$. Their composition is
an order two isometry that fixes a circle which does not intersect the
knot. Figure \ref{Vsymmetries} below shows the axis of the symmetries
$g_1$ and $g_2$. The axis of the order two symmetry
$g_1g_2$ is perpendicular to the page.

By the solution to the Smith conjecture, the fixed point set of
$g_1g_2$ in $S^3$ is unknotted circle which does not intersect the
knot.  The quotient of $S^3$ by  $g_1g_2$ is again $S^3$,
and the image of the two-bridge knot is another knot in $S^3$.  We
claim this knot is the unknot. Indeed, $g_1g_2$ is a symmetry whose
fixed point set is disjoint from the knot and which takes one bridge
to the other bridge.  The fundamental group of the two-bridge knot
complement is generated by two elements $x_1$ and $x_2$.  
  Fix a base point $b$ on the fixed point set of $g_1g_2$ in $B_1$.
  Then $x_1$ and $x_2$ can be represented by two curves which start at
  $b$ and encircle one of the two bridges in $B_1$. Now $g_1g_2$
acts on $\pi_1(S^3\setminus K)$ by setting $x_1= x_2$. Consider the
orbifold fundamental group  $\pi_1(P)$ where $P$ is $S^3 -
  K$ modulo the action of $g_1g_2$.  $\pi_1(P)$ is obtained by
adjoining an element $\omega$ to $\pi_1(S^3 - K)$ and adding the
relations $\omega^2 = 1$ and $\omega x_1 \omega = x_2$.  The
fundamental group of  $P$ is the quotient of  $\pi_1(P)$  by the normal closure of $\omega$.  Thus the
fundamental group of the underlying space of $P$ is
generated by one element.  As above, this is the
complement of a knot in $S^3$. Hence the knot is the unknot and the
underlying space  of $P$ is a solid torus. 
  This proves the claim.  The image of the singular set is also an
unknot, and it wraps around the image of the knot.

\begin{figure}[htb]
\label{orbifold}
\begin{center}
\input{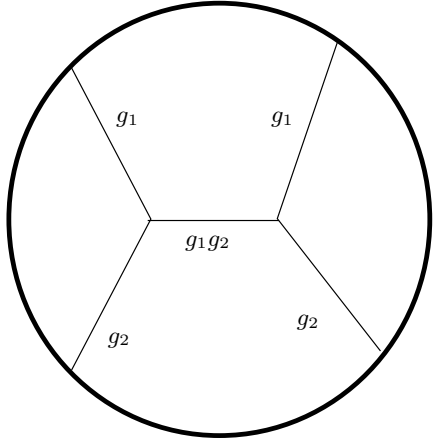}
\end{center}
\caption{$O$, the two-bridge knot complement modulo the action of $V$. The arcs of the singular set are order 2.}

\end{figure}

The fixed point sets of $g_1$, $g_2$ and $g_1g_2$ intersect in one
point in each of the $B_i$. The orbifold $O= (S^3 \setminus K)/ \V $
has the interior of a ball as its underlying space, and its cusp is a
copy of $S^2(2,2,2,2)$ (Figure \ref{orbifold} is a schematic picture).
The interior singular set of $O$ is a graph with two valence three
vertices where all arcs have order two singularities.  In Figure
\ref{orbifold}, the images of the fixed point sets of $g_1$, $g_2$ and
$g_1g_2$ in $O$ are labeled by their respective group elements.  We
will refer to the image in $O$ of the fixed point set of $g_1g_2$ by
$a_{g_1g_2}$.  Note that $a_{g_1g_2}$ will in general wind around the
arcs which meet the boundary.  However, the arcs which meet the
boundary are unknotted since unbranching along these arcs yields
 $P$, the quotient of the two-bridge knot complement by
$g_1g_2$, which has underlying space is a solid torus.\\[\baselineskip]
\noindent{\bf Proof of (i):}~In the case that $\Isom^+(S^3\setminus K) = V$, 
$O=Q$. A presentation of the orbifold fundamental
group $\Gamma$ of $Q$ can be obtained by removing the singular set,
taking a Wirtinger presentation of the complement, and then setting
$\gamma^2 = 1$ for each generator.  In this case, $N(K)$ is clearly
generated by elements of order 2.

In the case when  $\Isom^+(S^3\setminus K) =D_4$, $Q = {\bf H}^3/ N(K)$ is an orbifold which is double covered by
$O$. Call the generator of the covering group $\tau$.  By Theorem
\ref{symms}, no symmetry acts freely on $S^3 \setminus K$. Hence
$\tau$ has fixed points acting on $O$. Denote the fixed point set of
$\tau$ in $O$ by $b_\tau$.  Since $O$ is topologically a ball,
$b_\tau$ is a properly embedded arc, which is unknotted by the
solution to the Smith Conjecture.

$b_\tau$ meets the cusp of $O$ in two points.   Suppose one or both of these points coincide with the singular set of the cusp. Then there is an
element of order 4 in $N(K)$ that fixes a point on the cusp.  However,
any isometry of a knot complement takes a longitude to a longitude.
Hence if an isometry fixes a point on the cusp it can have order at
most 2.   Thus the fixed point set of $b_{\tau}$ is disjoint from the singular set of $O$ on the cusp at infinity. 

Since $\tau$ takes the singular set of $O$, a finite tree, to itself,
there is at least one fixed point in the interior of the singular set
of $O$.  We claim that there is exactly one, in the middle of
$a_{g_1g_2}$.  Since $\tau$ takes the singular set to itself, and
$b_\tau$ does not intersect the singular set in the cusp, the only
possible arc of intersection is $a_{g_1g_2}$.  In this case $b_\tau$
is an arc from the cusp to the cusp which strictly contains
$a_{g_1g_2}$.  The pre-image of $b_\tau$ in $S^3 \setminus K$ is the
fixed point set of some isometry $\tilde \tau$ of $S^3 \setminus K$.
This is either a circle or two arcs meeting the cusp. However, the
pre-image of $a_{g_1g_2}$ in $S^3 \setminus K$ is a circle, and a
circle or two arcs cannot properly contain a circle.  Therefore
$b_\tau$, the fixed point set of $\tau$, intersects the interior of
the singular set of $O$ in points.  Since $\tau$ takes the singular
set to itself, any such point must be in the center of the arc
$a_{g_1g_2}$, and there can only be one. This proves the claim. It
follows that $Q = O/ <\tau> = {\bf H}^3/N(K)$ is an orbifold with a
$S^2(2,2,2,2)$ cusp, and singular set consisting of a graph with order
2 arcs.  Combinatorially, the graph of the singular set of $Q$ is a
$H$, as is the graph of the singular set of $O$ in Figure
\ref{orbifold}.  Again, the two arcs of the singular set which meet
the boundary are unknotted.  One of the arcs is the image of $b_\tau$
in $Q$ which is unknotted.  The other arc is one unknotted arc of the
singular set of $O$ identified to another.  As before, using the
Wirtinger presentation of the ball with the singular set removed, we
see that $N(K) = \pi_1(Q)$
can be generated by elements with order 2. This proves (i).\\[\baselineskip]
\noindent{\bf Proof of (ii):}~As proved above, the cusp of $Q$ is a copy of
$S^2(2,2,2,2)$.  The only two-fold orbifold cover of this cusp which
is a torus is obtained by unbranching along the four cone points.
This is the kernel of a homomorphism
$\phi:\pi_1^{\orb}(S^2(2,2,2,2))\rightarrow
\mathbb{Z}/2\mathbb{Z}=\{\pm 1\}$ which sends each element of order
two to -1.  To extend this cover to $Q$, we must unbranch along the
curves of the singular set that meet the boundary.  The singular set
is combinatorially an H, and the underlying space is topologically a
ball.  By using a Wirtinger presentation for the orbifold group, it
follows easily that there can be no other unbranching.  Therefore
there is a unique two-fold cover of $Q$ with a torus cusp. We denote
this orbifold by $Q'$.  Since it was shown in the proof of (i) that
the unbranching arcs are unknotted in $Q$, it follows that $Q'$ has
underlying space a solid torus.  The singular set of the cover is the
pre-image of part of the singular set of $Q$, and so is of order 2.
Let $\pi_1^{\orb}(Q')$ denote the orbifold fundamental group of $Q'$,
and let $\psi:\pi_1^{\orb}(Q') \rightarrow \mathbb{Z}$ be the
homomorphism to the fundamental group of the underlying space of $Q'$.
If $\alpha$ is in the kernel of $\psi$, then $\alpha$ bounds an
immersed disk in the underlying space of $Q'$, and hence bounds an
immersed disk with $n$ order two singularities in $Q'$.  Represent
$\alpha$ by a 1-chain in some triangulation of $Q'$.  Then, after
subdividing, 2$\alpha$ bounds a $n$-punctured immersed sphere in $Q'$.
Thus $[\alpha] \in H_1(Q', \mathbb{Z})$ has order 2 when $\alpha$ is
in the kernel of $\psi$.  The map $\psi: \pi_1^{\orb}(Q') \rightarrow
\mathbb{Z}$ factors through the abelianization $H_1(Q', \mathbb{Z})
\rightarrow \mathbb{Z}$, so if $\alpha$ is not in the kernel of
$\psi$, $[\alpha]$ has infinite order. \qed

\section{Proof of Theorem \ref{main}} 

The proof of Theorem \ref{main} will follow immediately 
from the following two results.

\begin{theorem}
\label{nohiddensymms}
Let $K$ be a hyperbolic non-arithmetic 2-bridge knot.  Then $K$ has no
hidden symmetries.
\end{theorem}

Theorem \ref{nohiddensymms} together with Proposition \ref{nr}
shows that the minimal element in the orientable
commensurability class of $\Gamma_K$ is the group $N(K)=N(\Gamma_K)$. Hence $Q_K = {\bf H}^3/N(K)$
and the proof of Theorem \ref{main} is completed by:

\begin{theorem}
\label{onlyone}
Let $Q_K$ be as above. Then $Q_K$ is covered by exactly one knot complement in $S^3$.\end{theorem}

We defer the proof of Theorem \ref{nohiddensymms} and the remainder of
this section will be spent proving Theorem \ref{onlyone}. For
convenience we record the following result that will be used
subsequently in several places.  The non-trivial implication is
Theorem 3.4 (1) in \cite{GAW}.

\begin{theorem}  \label{cyclic}
  Let $\Sigma$, $\Sigma'$ be homotopy 3-spheres and $K \subset
  \Sigma$, $K' \subset \Sigma'$ non-trivial knots.  Then $\Sigma
  \setminus K$ is covered by $\Sigma' \setminus K'$ if and only if
  $\Sigma \setminus K$ admits a non-trivial cyclic surgery.
\end{theorem}

Of particular relevance to us is that in \cite{Ta} it is shown that a
hyperbolic 2-bridge knot complement does not admit a non-trivial surgery
with cyclic fundamental group.  Therefore, we have.

\begin{corollary}  \label{Taka} If $S^3 \setminus K$ is a hyperbolic 2-bridge knot complement,  it is not covered by the complement of any other knot in a homotopy 3-sphere. 
\end{corollary}

\noindent {\bf Proof of Theorem \ref{onlyone}}\\[\baselineskip]
Throughout the proof of Theorem \ref {onlyone}, $\Gamma_K$ will denote
a faithful discrete p-rep in standard form. Thus, the meridians $x_1$
and $x_2$ have images $\begin{pmatrix} 1 & 1 \\ 0 & 1 \end{pmatrix}$
and $\begin{pmatrix} 1 & 0 \\ y & 1 \end{pmatrix}$ respectively.

Assume $K'$ is a knot in $S^3$ such that $S^3\setminus K'$ is
commensurable with $S^3\setminus K$, and let $S^3\setminus K' = {\bf
  H}^3/\Gamma_{K'}$.  We can assume that $\Gamma_{K'}$ has been
conjugated to be a subgroup of $N(K)$.  Let $\Gamma = <\Gamma_K,\Gamma_{K'}>$
be the subgroup
of $N(K)$ that is generated by $\Gamma_K$ and $\Gamma_{K'}$, and let
$\Delta = \Gamma_K \cap \Gamma_{K'}$.  We have the following lattice of
subgroups of $\Gamma$. Note that
Theorem \ref{nohiddensymms} and Proposition \ref{nr} show that all 
the inclusions shown are of normal subgroups.
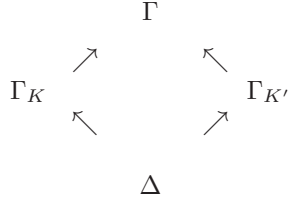
\eject
$$\Gamma$$

\centerline{$\nearrow$ \hskip .5 in $\nwarrow$}

\centerline{$ \Gamma_K$ \hskip 1 in $\Gamma_{K'}$}

\centerline{$\nwarrow$ \hskip .5 in $\nearrow$}

$$\Delta$$

Now $\Gamma_{K'}/\Delta \cong \Gamma/\Gamma_K $, which by Theorem
\ref{symms}, is a subgroup of $\V$ or $D_4$.  Since $\Gamma_{K'}$ is a
knot group, the quotient group $\Gamma_{K'}/\Delta$ cannot be $\V$ or
$D_4$ (since  both have non-cyclic abelianization).  Hence, the only
possibilities for $\Gamma_{K'} / \Delta \cong \Gamma/\Gamma_K$ are
cyclic of order $1$, $2$ or $4$.
Thus, since  ${\bf H}^3/\Delta$ is a cyclic cover of 
$S^3\setminus K'$, we deduce that ${\bf H}^3/\Delta$ has 1 torus cusp.

We claim that $\Gamma_K/\Delta$ is also cyclic.  To see this, we have
that ${\bf H}^3/\Delta \rightarrow S^3\setminus K$ is a regular cover
by a 1-cusped manifold.  Hence the peripheral subgroup of $\Gamma_K$
surjects the covering group. Therefore, the covering group is abelian,
and hence cyclic since $K$ is a knot.  Hence all of the inclusions indicated
in the above diagram have cyclic quotients.

Note first that if $\Gamma/\Gamma_K$ is the trivial group then $\Gamma_{K'}$
is a subgroup of $\Gamma_K$ which contradicts Corollary \ref{Taka}.
Thus we assume henceforth that $\Gamma\neq \Gamma_K$.

\begin{lemma} 
\label{toruscusp}
${\bf H}^3/\Gamma$ is not a manifold but has a torus 
cusp. \end{lemma}
\noindent {\bf Proof:}~Since $\Gamma\neq \Gamma_K$, 
the last sentence in the statement of Theorem \ref{symms} immediately
implies that ${\bf H}^3/\Gamma$ is an orbifold. Now suppose that ${\bf
  H}^3/\Gamma$ does not have a torus cusp.  Since it is the quotient
of a torus, and not rigid (by Theorem \ref{nohiddensymms} 
and Proposition \ref{nr}), it must be a  $S^2(2,2,2,2)$.
This implies that the
peripheral subgroup of $\Gamma$ fixing $\infty$ is
$$<\begin{pmatrix} 1 & 1 \\ 0 & 1 \end{pmatrix},\begin{pmatrix} 1 & r
  \\ 0 & 1\end{pmatrix}, \begin{pmatrix} i & a \\ 0 & -i
\end{pmatrix}>,$$
for some  numbers $r$ and $a$, and where
the third generator is the hyperelliptic involution of the torus that
takes a peripheral element to its inverse.

As remarked in \S 2.2, since $\Gamma$ is generated by parabolic
elements, we have that $\mathbb{ Q}(\tr \Gamma) = k\Gamma =
\mathbb{Q}(y)$.  We claim that this is a contradiction.

To that end, consider the following products of elements in $\Gamma$;
$$\begin{pmatrix} i & a \\ 0 & -i \end{pmatrix} \begin{pmatrix} 1 & 0 \\ y &1 \end{pmatrix}~\hbox{and}~ \begin{pmatrix} i & a \\ 0 & -i \end{pmatrix} \begin{pmatrix} 1 & 0 \\ y &1 \end{pmatrix}  \begin{pmatrix} 1 & 1 \\ 0 & 1 \end{pmatrix}.$$

The trace of the first product shows that $a\in \mathbb{Q}(y)$. The trace
of the second product then shows that $i\in \mathbb{Q}(y)$, in contradiction to Proposition
\ref{no2}. This completes the proof of the lemma. \qed

\begin{lemma} 
\label{equal}
$\Gamma/\Gamma_K \cong \Gamma/ \Gamma_{K'}$.
\end{lemma} 

\noindent {\bf Proof:}~The coverings 
$p_1:S^3\setminus K={\bf H}^3/\Gamma_K \rightarrow {\bf H}^3/\Gamma$
and $p_2:S^3\setminus K'={\bf H}^3/\Gamma_{K'} \rightarrow {\bf
  H}^3/\Gamma$ are cyclic coverings of an orbifold with a torus cusp
by a knot complement.  We denote by $\tau_1$, $\tau_2$ be the
generators of the cyclic group of covering transformations of the
covers $p_1$ and $p_2$ respectively. The lemma will follow on
establishing that each of $p_1$ and $p_2$ restricted to the cusp tori
$T_K$ and $T_{K'}$ of $S^3\setminus K$ and $S^3\setminus K'$
respectively have the same degree. To that end, first observe that the
preferred longitude $\ell = \begin{pmatrix} 1 & 2g_0 \\ 0 & 1
\end{pmatrix}$ of $S^3\setminus K$ is also a longitude for
$S^3\setminus K'$.  To see this, we have already noted that both of
these knot complements are cyclically covered by the 1-cusped manifold
${\bf H}^3/ \Delta$.  Since they are knot complements, their longitudes,  and the non-separating surfaces that these longitudes bound, both 
lift to this cyclic cover.    Since there is only one homology class in the cusp of ${\bf H}^3/ \Delta$ which bounds a non-separating surface, it follows that $\ell$ is also a
longitude of $S^3\setminus K'$.
By Lemma \ref{toruscusp} the orbifold ${\bf H}^3/
\Gamma$ has one torus cusp which we denote by $T_\Gamma$.  
Standard arguments in the orbifold setting (see the proof of
Theorem 11 of \cite{Du} for example)
provide a non-separating 2-orbifold in ${\bf H}^3/\Gamma$ which  is 
bounded by a simple closed curve in $T_\Gamma$. We denote this class 
(which is unique) by $\ell_0$. 
The pre-image of $[\ell_0 ]\in H_1(T_\Gamma; \mathbb{Z})$ must be $[\ell]
\in H_1(T_i; \mathbb{Z})$ in both coverings. We deduce from these remarks
that the covering degrees of $p_1$ and $p_2$ restricted to $\ell_0$ are
the same.

The proof will be completed by establishing that $x_1$ and $x_1'$,  the meridians of $\Gamma_K$ and $\Gamma_K'$, respectively,  are 
primitive elements of $\Gamma$.   Then the covering degree will be the degree restricted to $\ell_0$. To prove primitivity we argue as follows.
Note that $\Gamma_K$ has algebraic 
integer entries.  We first claim that $\Gamma$, and hence $\Gamma_{K'}$,
also has algebraic integer entries. Since $\Gamma$ has one torus cusp we can
choose left coset representatives of $\Gamma_K$ in $\Gamma$ that are
parabolic of the form $\begin{pmatrix} 1 & a \\ 0 & 1
\end{pmatrix}$.  Since the property of having algebraic integer traces is a
commensurability invariant, $\text{tr}\begin{pmatrix} 1 & 0 \\ y & 1 \end{pmatrix}
\begin{pmatrix} 1 & a \\ 0 & 1 \end{pmatrix} = 2 + ay$
is also an algebraic integer.  Since $y$ is a unit ($\Lambda_K(y)$
is  monic with constant term $1$) $a$ is also an
algebraic integer. Since this holds for any left coset
representative, this proves the claim.

Now $x_1 = \begin{pmatrix} 1 & 1 \\ 0 & 1 \end{pmatrix}$, and we shall
assume that 
$x_1' = \begin{pmatrix}
  1 & r \\ 0 & 1 \end{pmatrix}$. Note that $r$ must be a unit, for
if not, since we have shown all entries of $\Gamma_{K'}$ are
integral, we can find a prime ideal $\cal P$ 
dividing $<r>$, and reducing the entries modulo $\cal P$ sends all of 
$\Gamma_{K'}$ to the identity (since it is normally generated by $x_1'$). However,
this impossible since there is an element in $\Gamma_{K'}$ conjugating
$x_1'$ to a meridian $x_2'$ fixing 0. Such a conjugating
matrix has zero as the $(1,1)$ entry.

If $x_1'$ is not primitive, there is an element $\begin{pmatrix} 1 & t
  \\ 0 & 1 \end{pmatrix} \in \Gamma$ such that
$$\begin{pmatrix} 1 & t \\ 0 & 1 \end{pmatrix}^n = \begin{pmatrix} 1 &
  r \\ 0 & 1 \end{pmatrix},$$
for some integer $n$ . We claim $n=1$.
To see this, note that  the elements of $\Gamma$ have algebraic integer entries 
by the argument above.  Therefore $t=r/n$ is an algebraic integer.
Since $r$ is a unit, $t=\pm r$, and $x_1'$ is not a proper power of any element
in $\Gamma$.  It is clear that the same
argument holds for $x_1$.  Therefore, both $x_1$
and $x_1'$ are primitive in $\Gamma$ as required. \qed\\[\baselineskip]
%Also, the action of $\tau_1$ and $\tau_2$ on the cusp tori $T_K$
%and $T_{K'}$ induces
%an action on $H_1(T_K; \mathbb{Z})$ and $H_1(T_{K'}; \mathbb{Z})$ that
%takes the homology class of the longitude $\ell$ to itself in each
%case.  $\tau_1$ and $\tau_2$ must also take the class of the meridian
%$x_1$ of $K$ and $x_1'$ of $K'$ to a class which is dual to $\ell$.
%Since $\tau_1$ and $\tau_2$ have finite order, it follows that
%their action takes $x_1$ to $x_1$ and $x_1'$ to $x_1'$ respectively. 
From the previous discussion,  we are assuming that $\Gamma/\Gamma_K =
\Gamma/\Gamma_{K'}$ is a cyclic group of order 2 or 4. 

In the case when $\Gamma/\Gamma_K = \Gamma/\Gamma_{K'} =
\mathbb{Z}/4 \mathbb{Z}$ this follows easily from Lemma \ref{oneZ4}.   Then since ${\bf H}^3/\Gamma$ has a torus cusp, 
$ {\bf H}^3/\Gamma$ is  the unique two-fold cover  $Q'$ of ${\bf H}^3 /
N(\Gamma_K)$ with a torus cusp, and Lemma \ref{oneZ4}, shows
that ${\bf H}^3 / \Gamma$ has exactly one 4-fold cyclic cover. 
Therefore $\Gamma_K = \Gamma_{K'}$.  \\[\baselineskip]
\noindent We now assume that $\Gamma/\Gamma_K = \Gamma/\Gamma_{K'} = 
\mathbb{Z}/2 \mathbb{Z}$.
In this case all the cyclic quotients arising
from the lattice of subgroups at the beginning of the proof of Theorem
\ref{onlyone} are order 2. In particular, since $K$ is 2-bridge, ${\bf
  H}^3/ \Delta$ is the complement of a knot in a lens space.

Recall the proof of Lemma \ref{equal} shows that $x_1$ and $x_1'$
are primitive elements in $\Gamma$.  Furthermore, $x_1 \notin
\Gamma_K'$ and $x_1' \notin \Gamma_K$.
For if so, then $x_1 \in \Gamma_K \cap \Gamma_{K'} = \Delta$ which is
false. If $x_1'\in \Gamma_K$ then $x_1'\in \Delta$ and normality implies
that $\Delta=\Gamma_{K'}$ which is false as $[\Gamma_{K'}:\Delta]=2$.
%and that the cyclic covers in question
%are determined by the action on $\ell$. From the description of $\ell$
%given in \S 2.3, it follows that a longitude for the cusp torus of
%${\bf H}^3/\Gamma$ is $\ell_0 = \begin{pmatrix} 1 & g_0 \\ 0 &
%1 \end{pmatrix}$ (which is also primitive in $\Gamma$).  
%Since $[\Gamma:\Gamma_K]=2$ and $\ell_0\notin\Gamma_K$ it follows
%that $\Gamma=<\Gamma_K,\ell_0>$.
Consider the normal closure of $x_1^2$ in
$\Gamma$. We denote this by $<x_1^2>_\Gamma$. We claim that
$<x_1^2>_{\Gamma} = <x_1^2>_{\Gamma_K}$.
The inclusion $<x_1^2>_{\Gamma_K}\subset <x_1^2>_{\Gamma}$ is clear. For
the reverse inclusion, 
we can choose $x_1'$ to be a left coset representative 
for $\Gamma_K$ in
$\Gamma$ and it follows that $\Gamma=<\Gamma_K, x_1'>$.  Since
$x_1'$ commutes with $x_1^2$, we deduce that
$<x_1^2>_\Gamma \subset <x_1^2>_{\Gamma_K}$.

The exact same argument with $\Gamma_K$ and $x_1'$ replaced by
$\Gamma_K'$ and $x_1$ shows that $<x_1^2>_\Gamma =
<x_1^2>_{\Gamma_K'}$.  Similarly, we can take $x_1$ to be a left coset
representative of $\Delta$ in $\Gamma_K$, which also shows that
$ <x_1^2>_\Delta
=<x_1^2>_{\Gamma_K}$. Thus we conclude that 
$<x_1^2>_{\Gamma_K} = <x_1^2>_{\Gamma_K'} = <x_1^2>_\Delta$.  
For convenience we denote $<x_1^2> _{\Gamma}$ by
$\cal N$.  

Now the group $\Gamma_K/{\cal N}$ is the orbifold fundamental group of
the orbifold obtained by the $(2,0)$ Dehn filling on the 2-bridge knot
$K$. Now the double branched cover of $K$ is
a lens space $L(p/q)$ whose fundamental group is the cyclic group of
(odd) order $p$. Since 
this double cover is
obtained by first performing $(2,0)$ orbifold Dehn surgery on $K$, and
then passing to the index $2$ cover which is a manifold, we deduce that
$\Gamma_K/{\cal N}$ is a dihedral group of order $2p$ for
some odd integer $p$.  

Note that $x_1^2$ is a primitive element of $\Gamma_{K'}$. For if not, then 
$x_1 \in \Gamma_{K'}$, but as noted above 
$x_1\notin \Delta$. Therefore, $\Gamma_K'/{\cal N}$ and
$\Delta/{\cal N}$ are the fundamental groups of the
manifolds obtained by Dehn filling the primitive curve $x_1^2$ in $S^3\setminus
K'$ and ${\bf H}^3/\Delta$ respectively. In particular, this latter
group is a cyclic group of order $p$ (denoted $C_p$). We denote the former
group by $G'$.
Hence we have the following diagram of groups: \\[\baselineskip]

 \centerline{$G=\Gamma/{\cal N}$}

\centerline{$\nearrow$ \hskip .5 in $\nwarrow$}

\centerline{$  D_p$ \hskip 1 in $G'$}

\centerline{$\nwarrow$ \hskip .5 in $\nearrow$}

\centerline{$C_p$}

\vskip .3 in

The group $G$ is a 2-fold extension of $D_p$, and so has order $4p$.  
The group $G'$ is a finite group of order $2p$
arising as the fundamental group of
a closed orientable 3-manifold. We claim that the finite group $G'$ is cyclic.

This will complete the proof of Theorem \ref{onlyone}. For if $G'$ is
cyclic, then $S^3\setminus K'$ has a
$2p$-fold cyclic cover $\Sigma\setminus K''$, where $\Sigma$ is a
homotopy 3-sphere.  Note that since $\Sigma\setminus K''$ is the
$2p$-fold cylic cover of $S^3\setminus K'$, it must cover the 2-fold
cover of $S^3\setminus K'$, namely ${\bf H}^3/\Delta$ . Hence
$\Sigma\setminus K''$ covers $S^3\setminus K$ which contradicts
Corollary \ref{Taka}.

To  establish that $G'$ is cyclic, by \cite{Mi} there are a limited number of types of
non-cyclic finite groups that can be the fundamental group of a closed
orientable 3-manifold. We list these below using the notation of
\cite{BZ}:
\begin{itemize}
\item {\em Even $D$-type: $\lbrace x,y: x^2 = (xy)^2 =y^n \rbrace
    \times \mathbb{Z}/j\mathbb{Z}$, with $j \geq 1,\ \ n\geq 2$ with
    $n$ even.  The abelianization is $\mathbb{Z}/2\mathbb{Z} \oplus
    \mathbb{Z}/2j\mathbb{Z}$.}
\item {\em Odd $D$-type: $\lbrace x,y: x^{2^k} =1,y^{2l +1} = 1,
    xyx^{-1} = y^{-1} \rbrace \times \mathbb{Z}/j\mathbb{Z}$, with $j \geq 1$,
$(2(2l+1),j) =1$ $k \geq 2$.  The
    abelianization is $\mathbb{Z}/2^k j\mathbb{Z}$.}
\item {\em $T$-type: $\lbrace x,y,z: x^2 = (xy)^2 = y^2, z^{3^k} =1,
    zxz^{-1} =y, zyz^{-1} = xy \rbrace \times \mathbb{Z}/j\mathbb{Z},
    (6,j) =1$.  The abelianization is $\mathbb{Z}/3^k j\mathbb{Z}$.}
\item {\em $O$-type: $\lbrace x,y: x^2 = (xy)^3 = y^4, x^4 = 1 \rbrace
   \times \mathbb{Z}/j\mathbb{Z} , (6,j) = 1$.  The abelianization is
    $\mathbb{Z}/2j\mathbb{Z}$.}
\item {\em $I$-type: $\lbrace x,y: x^2 = (xy)^4= y^5, x^4 =1 \rbrace
    \times \mathbb{Z}/j\mathbb{Z}, (30,j)=1$.  The abelianization is
    $\mathbb{Z}/j\mathbb{Z}$.}
\item {\em $Q$ type: $\lbrace x,y,z: x^2 = (xy)^2 = y^{2n}, z^{kl} =
    1, xzx^{-1} = z^r, yzy^{-1} = z^{-1} \rbrace \times
    \mathbb{Z}/j\mathbb{Z}$, $n, k, l, j$ relatively prime odd
    positive integers, $r \equiv -1 \mod k , r \equiv 1 \mod l$.  The
    abelianization is $\mathbb{Z}/2\mathbb{Z} \oplus
    \mathbb{Z}/2j\mathbb{Z}$.}

\end{itemize}

Now $|G'| = 2p$ where $p$ is odd. This allows us to immediately rule
out both $D$-types, $O$-type, $I$-type and $Q$-type, since either the
fundamental group or its abelianization is clearly divisible by 4.
Since $G'$ has order $2p$ and contains a cyclic subgroup of order $p$, $G'$
surjects $\mathbb{Z}/2 \mathbb{Z}$. This precludes $G'$ being
of T-type, since the abelianization of a group of $T$-type
is odd. This proves that $G'$ is a cyclic group as claimed.\qed

\section{Two-bridge knots and hidden symmetries}

In this section we prove Theorem \ref{nohiddensymms}.  This will be
done in \S 4.1 and 4.2.  We begin with some preliminary discussion.

\subsection{}  
$K$ is a hyperbolic 2-bridge knot different from the figure-eight knot,
and as above,  $\Gamma_K$ the faithful discrete p-rep of
$\pi_1(S^3\setminus K)$ in standard form; ie given by:

$$\rho(x_1) = \begin{pmatrix} 1 & 1 \\ 0 & 1 \end{pmatrix}~\hbox{and}~\rho(x_2) = 
\begin{pmatrix} 1 & 0 \\ y & 1 \end{pmatrix}.$$

As above the invariant trace-field is $\mathbb{Q}(y)$, and we let the ring of integers of $\mathbb{Q}(y)$, be denoted by $R_y$.
Assuming that $K$ has hidden symmetries,  Corollary \ref{noi} shows that
the orientable commensurator orbifold $Q_K$ has a $S^2(2,3,6)$ rigid
cusp. Hence  $\mathbb{Z}[\omega] \subset R_y$. 

Let $\wp \subset R_y$ be a prime ideal that divides the principal
ideal $2R_y \subset R_y$, and let ${\bf F} = R_y/\wp$.  This is a
finite field of order $2^s$ for some integer $s\geq 1$. In fact, since
$\mathbb{Z}[\omega]\subset R_y$ and $2$ is inert in
$\mathbb{Z}[\omega]$, it follows that $|{\bf F}|=4^s$ for some integer
$s\geq 1$.  Let $\phi:\PSL(2,R_y) \rightarrow \PSL(2,{\bf F})$ be the reduction
homomorphism. The key result that is needed to prove Theorem
\ref{nohiddensymms} is the following.

\begin{theorem} 
\label{imagecontrol}
The image of $\Gamma_K$ under $\phi$ is a dihedral group of order 6 or 10. 
\end{theorem}

Deferring the proof of Theorem \ref{imagecontrol} until \S 4.2, we complete the proof of Theorem \ref{nohiddensymms}.

As in \S 2.3 we let $\Lambda_0(y)$ denote the factor of the p-rep
polynomial $\Lambda_K(y)$ that corresponds to the complete hyperbolic
structure.  Let $\bar \Lambda_0(y)$ be the reduction of $\Lambda_0(y)$
modulo 2. By Theorem \ref{norepeated}, $\bar \Lambda_0(y)$ has no
repeated factors, and each factor corresponds to a dihedral
representation of $\Gamma_K$. By Theorem \ref{imagecontrol} the image
of these representations are dihedral of order $6$ or $10$. Now from
\cite{Ri1} Proposition 4, we deduce that the corresponding factors of
$\bar \Lambda_0(y)$ have degrees 1 and 2 respectively.  Since we are
working modulo 2, the only possible irreducible polynomials are $x$ in
the case of $D_3$ and $x^2 + x + 1$ in the case of $D_5$.

 By \ref{norepeated} and \ref{numbertheory}(i) 2 does not divide the discriminant  $\Delta(y)$ of $\Lambda_0(y)$. Hence 2 does not divide $\Delta(y)d_k^{-1}$, where $d_k$ is the discriminant of $\mathbb{Q}(y)$. Therefore Theorem \ref{numbertheory}(ii) 
shows that the above factors determine the
decomposition of the principal ideal $2R_y$ into prime (ideal) power
factors.  Any such $\wp$ as above is one of these factors.  Since there are no repeated factors, we deduce that the
decomposition of $2R_y$ is into at most the product of two prime
ideals of $R_y$.  Therefore, $\bar\Lambda_0(y)$ has degree at most
$3$. Hence $\Lambda_0(y)$ has degree at most $3$ (recall
$\Lambda_0(y)$ is a monic polynomial so the leading coefficient is
never zero modulo 2).

The degree cannot be 1, since $y$ is a non-real root of the p-rep polynomial.
If the degree is 2, then since $\mathbb{Q}(y)$ contains $\mathbb{Q}(\sqrt{-3})$
we have that $\mathbb{Q}(y) =\mathbb{Q}(\sqrt{-3})$ and so $\Gamma_K$ has traces in $\mathbb{Z}[\omega]$. It follows that $\Gamma_K$ is arithmetic (cf. \cite{Re1}) which is
false.  Finally, the degree cannot be 3, because $\mathbb{Q}(\sqrt{-3}) \subset \mathbb{Q}(y)$.
This contradiction completes the proof.\qed

\subsection{Proof of Theorem \ref{imagecontrol}}

Throughout this subsection we let $C(K)=C^+(\Gamma_K)$
which is assumed to contain $\Gamma_K$. Since $\Gamma_K$ contains a
peripheral subgroup fixing $\infty$, $C(K)$ has a peripheral
subgroup fixing $\infty$ and we denote this by $B$ (so $B$ is
isomorphic to the orbifold group of $S^2(2,3,6)$). We begin with some preliminary lemmas.

\begin{lemma} \label{form} The orbifold group $B<C(K)$ has the form
$$<\begin{pmatrix} 1 & 1 \\ 0 & 1 \end{pmatrix},
                    \begin{pmatrix} 1 & \omega \\ 0 & 1\end{pmatrix}, 
\begin{pmatrix} i\omega & 0 \\ 0 & -i\overline{\omega} \end{pmatrix}>.$$
\end{lemma}

\noindent{\bf Proof:}~$B$ fixes infinity 
and so an element of order $6$ in $B$ has the form $e=\begin{pmatrix}
  i\omega & t \\ 0 & -i\overline{\omega}
\end{pmatrix}$ for some  number $t$. We claim that we can
arrange that $t$ can be taken to be $0$.

To see this we argue as follows.  By \cite{NR} Theorem 2.2, there is a
normal subgroup $L$ of $C(K)$ (with quotient $(\mathbb{Z} /
2\mathbb{Z})^{\alpha}$) such that all elements of $C(K)$ whose trace
lies in $\mathbb{Q}(y)\setminus\{0\}$ are elements of $L$.  Recall
from the proof of Proposition \ref{no2} the element $\mu\in N(K) <
C(K)$. Hence the element $\delta=e\mu=\begin{pmatrix} -\omega & -ti \\ 
  0 & -\overline{\omega} \end{pmatrix}$ is an element of order $3$
which lies in $L$. Now $L$ also contains $\Gamma_K$ and so,
$$\tr(\rho(x_2)\delta) = -\omega-\overline{\omega}-yti \in \mathbb{Q}(y).$$
Furthermore, since  traces of elements in $\Gamma_K$ are algebraic integers,
all traces of elements in $C(K)$ are algebraic integers by commensurability. Letting $R_y$
denote the ring of integers in $\mathbb{Q}(y)$ we deduce
that $yti\in R_y$. Since $y$ is a unit we deduce that $ti\in R_y$.

Now $e^3 = \begin{pmatrix} -i & 2t \\ 0 & i \end{pmatrix}$ is an
element of order $2$ in $B$, and the product $\mu e^3 =
\begin{pmatrix} 1 & 2ti \\ 0 & 1 \end{pmatrix}$ is a parabolic element
in $B$.  Since the cusp field is $\mathbb{Q}(\sqrt{-3})$ it follows
that $ti\in \mathbb{Q}(\sqrt{-3})$.  Furthermore from above $ti$ is an
algebraic integer, and so $ti\in \mathbb{Z}[\omega]$.  Hence the element
$\delta$ above has coefficients in $\mathbb{Z}[\omega]$.  Let $x =
-ti\overline{\omega}\in \mathbb{Z}[\omega]$, and consider the product
$\delta\begin{pmatrix} 1 & x\\ 0 & 1 \end{pmatrix}$. This gives the
element $\begin{pmatrix} \omega & 0 \\ 0 & \overline{\omega}
\end{pmatrix}$. Taking the product with $\mu$
gives the desired element of order $6$.\qed

\begin{lemma}
\label{333cusp}
$C(K)$ contains a subgroup $C_0(K)$ of index $2$ such that $\Gamma_K$
is a subgroup of $C_0(K)$ and ${\bf H}^3/ C_0(K)$ has a $S^2(3,3,3)$
cusp. Furthermore, $C_0(K) <\PSL(2,R_y)$.\end{lemma}

\noindent{\bf Proof:}~As mentioned in the proof of 
Lemma \ref{form}, \cite{NR} Theorem 2.2 provides
a normal subgroup $L$ of $C(K)$ with quotient $(\mathbb{Z} /
2\mathbb{Z})^{\alpha}$ such that all elements of $C(K)$ whose trace
lies in $\mathbb{Q}(y)\setminus\{0\}$ are elements of $L$.  As noted,
$L$ contains $\Gamma_K$, and so ${\bf H}^3/L$ has one cusp. Since
${\bf H}^3/L$ has one cusp, $B$ must surject the covering group
$(\mathbb{Z} / 2\mathbb{Z})^{\alpha}$.  We claim that this forces
$\alpha = 1$.

To see this, note that the abelianization of the group $B$ is
$\mathbb{Z}/6\mathbb{Z}$. Hence, the image of $B$ under
the homomorphism $C(K)\rightarrow C(K)/L$ is cyclic, and 
so $\alpha\leq 1$. However, $\alpha\neq 0$ 
since an element of order $6$ in $B$ cannot lie in $L$. 
For if so, then $\sqrt{3}\in \mathbb{Q}(y)$.
Since $\mathbb{Q}(\sqrt{-3}) \subset \mathbb{Q}(y)$ 
(by Corollary \ref{noi}), it follows
that $\sqrt{3}\sqrt{-3}=3i\in \mathbb{Q}(y)$ 
which contradicts Proposition \ref{no2}.  
Given that the element $\begin{pmatrix} \omega & 0 \\ 0 & \overline{\omega} \end{pmatrix}$ constructed
in the proof of Lemma \ref{form} is an element of $L$, this argument also shows that the cusp of ${\bf H}^3/L$ is $S^2(3,3,3)$. 

Note that $L< \PSL(2,R_y)$. For since ${\bf H}^3/L$ has one cusp, a system of coset representatives
of $\Gamma_K$ in $L$ can be taken from the cusp subgroup $E$ of $L$.  By definition of $B$, the
subgroup $E$ is:

$$<\begin{pmatrix} 1 & 1 \\ 0 & 1 \end{pmatrix},
                    \begin{pmatrix} 1 & \omega \\ 0 & 1\end{pmatrix}, 
\begin{pmatrix} \omega & 0 \\ 0 & \overline{\omega} \end{pmatrix}>.$$

Thus all the coefficients lie in $\mathbb{ Z}[\omega]$. Now Corollary
\ref{noi} implies that $\mathbb{ Z}[\omega] \subset R_y$.
Hence all the coefficients of $L$ are elements of $R_y$. We now take $C_0(K)=L$.\qed\\[\baselineskip]
%Recall the group
%$C_0(K)<\PSL(2,R_y)$ (by Lemma \ref{333cusp}) contains $\Gamma_K$ and
%has a peripheral subgroup:
%$$E = <\begin{pmatrix} 1 & 1 \\ 0 & 1 \end{pmatrix},
%                    \begin{pmatrix} 1 & \omega \\ 0 & 1\end{pmatrix}, 
%\begin{pmatrix} \omega & 0 \\ 0 & \overline{\omega} \end{pmatrix}>.$$$  
Consider the reduction homomorphism $\phi$ restricted to the group
$C_0(K)$ (which is a subgroup of $\PSL(2,R_y)$ by Lemma
\ref{333cusp}). We continue to denote this by $\phi$ and let
$\Delta=\ker \phi$.  As in the proof of Lemma \ref{333cusp},
$E=<\begin{pmatrix} 1 & 1 \\ 0 & 1 \end{pmatrix},
                    \begin{pmatrix} 1 & \omega \\ 0 & 1\end{pmatrix}, 
\begin{pmatrix} \omega & 0 \\ 0 & \overline{\omega} \end{pmatrix}>$ denotes the cusp subgroup of the $C_0(K)$.
Note that $E \cap \Delta$ is torsion free, since the element
$\begin{pmatrix} \omega & 0 \\ 0 & \overline{\omega} \end{pmatrix}$ of
order $3$
injects under $\phi$. \\[\baselineskip]
We claim that $\phi(E)$ has order $12$. Indeed, since
$E<\PSL(2,\mathbb{ Z}[\omega])<\PSL(2,R_y)$ it suffices to consider
the image of $E$ under the reduction homomorphism restricted to
$\PSL(2,\mathbb{ Z}[\omega])$. By the definition of $E$ and $\phi$ it
is easily seen that $\phi(E)$ is an extension of $\V$ by $\mathbb{
  Z}/3\mathbb{ Z}$. This defines a subgroup of order $12$, which
proves the claim.

Now $\Gamma_K \cap \Delta$ is the kernel of $\phi$ restricted to
$\Gamma_K$. As discussed
in \S 2.3 this is a dihedral representation. 
We need to show that $\Gamma_K / \Gamma_K \cap \Delta$ 
is a dihedral group of order 6 or 10. 

\begin{lemma} \label{odd} $\Gamma_K / \Gamma_K \cap \Delta$ is a dihedral group of order $2m$
and $m$ is odd.\end{lemma} 

\noindent{\bf Proof:}~That $\Gamma_K \cap \Delta$ is dihedral follows from the sentence
before the lemma.  Since $K$ is a 2-bridge knot it has a 2-bridge normal
form as discussed in \S 2.3, and hence the double branched cover of $K$ is
a lens space $L(p/q)$ whose fundamental group is the cyclic group of
odd order $p$.

The double cover of a manifold branched over a knot $K$ can be
obtained by first performing $(2,0)$ orbifold Dehn surgery on $K$, and
then passing to the index $2$ cover which is a manifold.  In
particular meridians of $K$ are mapped to elements of order $2$, and
furthermore any quotient of $\pi_1(S^3\setminus K)$ in which the
meridians are mapped to elements of order $2$ is a quotient of the
orbifold group obtained above.

Hence, in the case at hand, the orbifold group is a dihedral group of
order $2p$
where $p>1$ is odd. The lemma
now follows from the discussion in the previous paragraph.\qed\\[\baselineskip]
We now complete the proof of Theorem \ref{imagecontrol}.
As before, if $\ell$ denotes the longitude for $x_1$
described in \S 2.3,
then $\rho(\ell)=\begin{pmatrix} 1 & g \\ 0 & 1 \end{pmatrix}$, where
$g = 2 g_0$ and $g_0\in R_y$.  Hence $\phi(\rho(\ell)) = 1$, and so
the image of the peripheral subgroup $<\rho(x_1),\rho(\ell)>$ under
$\phi$ is cyclic of order $2$. Hence the cover of $S^3\setminus K$
determined by $\Gamma_K\cap \Delta$ has $m$ cusps (and $m$ is odd by
Lemma \ref{odd}).

We now count the cusps of ${\bf H}^3/(\Gamma_K\cap \Delta)$ in a
different way.  Consider the following diagram of subgroups of $C_0(K)$:

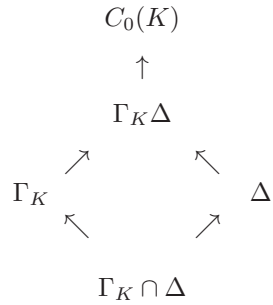
\medskip

%\vskip 1 in
$$C_0(K)$$
$$\uparrow$$
$$\Gamma_K \Delta$$

\centerline{$\nearrow$ \hskip .5 in $\nwarrow$}

\centerline{$  \Gamma_K$ \hskip 1 in $\Delta$}

\centerline{$\nwarrow$ \hskip .5 in $\nearrow$}

$$\Gamma_K \cap \Delta$$

Since $\Gamma_K \Delta$ contains $\Gamma_K$ and is contained in
$C_0(K)$, it has one cusp which is either a torus or $S^2(3,3,3)$.
Furthermore, $\Gamma_K\Delta / \Delta \cong \Gamma_K / (\Gamma_K \cap
\Delta) = D_{m}$.  We claim that the cusp of ${\bf
  H}^3/\Gamma_K\Delta$ is a torus. To see this we argue as follows.
Assume that the cusp is $S^2(3,3,3)$.  Denote the cusp subgroup of
$\Gamma_K\Delta$ fixing $\infty$ by $E'$.  Arguing as above for
$\phi(E)$ shows that $\phi(E')$ has order 12. Since 
$\phi(\Gamma_K\Delta)=D_{m}$, 12 must divide $2m$. However,
this contradicts Lemma \ref{odd} which shows $m$ is odd.

Therefore the cusp of $\Gamma_K \Delta$ is a torus, and we have
a covering space $S^3 \setminus K \rightarrow {\bf H}^3 / \Gamma_K
\Delta$. By \cite{Re1} Lemma 4 this is a regular abelian cover. By
Theorem \ref{symms} it follows that the covering group is cyclic of
order $1$,$2$ or $4$ or it is the group $\V$.

The image of $C_0(K)$ under $\phi$ is a subgroup of $\PSL(2, {\bf F})$, where 
as discussed in \S 4.1, $|{\bf F}| = 4^s=q$.  
 By \cite{Su} Theorem 6.25, subgroups of $\PSL(2, {\bf F})$ are as follows: 

\begin{enumerate}
\item $\PSL(2,{\bf F}')$ where ${\bf F}'$ is a subfield of ${\bf F}$
  of order $2^k$. Note that since the characteristic is $2$,
  $\PGL(2,{\bf F}')=\PSL(2,{\bf F}')$ which excludes one of the
  possibilities of \cite{Su} Theorem 6.25.
%by p. 414 of \cite{Su}.
\item
$A_4$ or $A_5$. Note that $S_4$ is ruled out by \cite{Su} Theorem 6.26(C).
\item
A subgroup $H$ of order $q(q-1)$ and its subgroups. A Sylow-2 subgroup $Q <H$ is elementary abelian and $H/Q$ is cyclic of order $q-1$. 
\item
A dihedral group of order $2(q \pm 1)$ or one of its subgroups.
\end{enumerate}

We handle each possibility in turn.\\[\baselineskip]
\noindent{\bf Case 1:}~Assume that $C_0(K)$ maps onto $\PSL(2, {\bf F}')$. 
From above we know that $[\Gamma_K\Delta:\Gamma_K]=1,2,4$. If
$[\Gamma_K\Delta:\Gamma_K]=1$ then $\Delta = \Gamma_K \cap \Delta$.
Now the number of cusps of ${\bf H}^3/ \Delta$ is given by
$|\PSL(2,{\bf F}')|/|\phi(E)|$, which is $2^k(2^{2k} -1) /12$ (by
our computation of $|\phi(E)|$ given above ).  From above, the number of cusps is the odd integer $m$. It
follows that $k = 2$, from which it follows that $m = 5$, thus
yielding the dihedral group $D_5$.

In the case that $\Gamma_K \Delta /\Gamma_K$ is cyclic of order 2 or 4
or $\V$, we have $\Delta/(\Gamma_K \cap \Delta)$ is cyclic of order 2
or 4 or $\V$. Thus we count that ${\bf H}^3/\Delta$ can have $m$,
$m/2$, or $m/4$ cusps.  Since $m$ is odd, $\Delta$ has $m$ cusps. We
now argue as above, and we deduce that
$k = 2$ and  $m = 5$ yielding the dihedral group $D_5$ again. \\[\baselineskip]
\noindent{\bf Case 2:}~Assume that $C_0(K)$ maps onto $A_4$ or $A_5$.
The only non-cyclic dihedral subgroup of $A_4$ is $\V$ and $\Gamma_K$ cannot map onto this
since it is a knot group. 
The only non-cyclic dihedral subgroups of $A_5$ are $D_{5}$ and $D_3$.  Hence we are done in this case.\\[\baselineskip]
\noindent{\bf Case 3:}~Suppose that the image of $C_0(K)$ is a subgroup of $H$ as given above.
$H$ contains the image of $\Gamma_K$. This is a dihedral subgroup
$D_m$ of order $2m$, where $m$ is odd.  Since $Q$ is normal in $H$, $Q
\cap D_m$ is normal in $D_m$, and the
only normal subgroup of $D_m$ is the cyclic subgroup of  order $m$.  Since $Q$ is a 2-group, 2 must divide $m$, but $m$ is odd. \\[\baselineskip]
\noindent{\bf Case 4:}~Suppose that the image of $C_0(K)$ is a subgroup of a dihedral group of order $2(q \pm 1)$. The order of the image of $E$ is 12. Hence 12 divides $2(4^s \pm 1)$, which is absurd. 
This completes the proof.\qed

\section{Commensurable knot complements}

\subsection{} In the non-hyperbolic case, infinitely many knot complements 
can easily occur in one commensurability class.  For example,
non-hyperbolic two-bridge knots (ie the 2-bridge torus knots) all have
commensurable complements (see \cite{NeumannCommen}). 

As mentioned in 
\S 1, there are examples of hyperbolic knots $K$ for which the commensurability
class of $S^3\setminus K$ contains more than one knot complement. 
We now describe the constructions that are known to us.\\[\baselineskip]
\noindent{\bf Lens space surgeries:}~The main source of examples
occurs when $K$ admits a lens space surgery. In this case $S^3\setminus
K$ has a cyclic cover that is a knot complement by Theorem \ref{cyclic}.
By \cite{CGLS}, there can be at most two non-trivial
cyclic surgeries.  Hence if a hyperbolic knot $K$ admits a lens space
surgery there can be up to 3 knot complements in the
commensurability class that arise from this construction. This holds for the $(-2,3,7)$-pretzel knot
complement (see \cite{Be}).  \\[\baselineskip]
\noindent{\bf Hidden symmetries:}~A pair of commensurable knot complements
that do not arise as above are the two dodecahedral knot
complements of \cite{AR}.  These are commensurable of the same
volume, and have hidden symmetries (see \cite{NR}).\\[\baselineskip]
If a hyperbolic knot $K$ has no hidden symmetries 
then there are only finitely 
knot complements in the commensurability class of $S^3\setminus K$
(see \cite{Re1} Theorem 5). However, in the presence of hidden
symmetries, it is unknown to the authors whether there are finitely
many knot complements in a commensurability class (even for the
dodecahedral knots).\\[\baselineskip]
\noindent{\bf Symmetries:}~The final construction we are aware of 
was described to us by W. Neumann (personal communication).
He has constructed an infinite family of pairs of knots $\{K_i,K_i'\}$ 
which have the following property. The complements $S^3\setminus K_i$
and $S^3\setminus K_i'$ have different volumes, and are
both regular covers of a common (genuine) orbifold. 

The simplest pair
of examples are the knots $9_{48}$ and $12n642$, where the volume ratio
is $4:3$. 
These examples do not arise in connection with lens space surgeries 
on a knot. Indeed, neither of these knots admits a lens space
surgery by \cite{WZ}. 

Since the trace-field of the knot complements associated to the knots
$9_{48}$ and $12n642$ is cubic, there are no hidden symmetries (recall
\S2.2), and so the remarks above show that there are finitely many
knot complements in this commensurability class (presumably 2).

\subsection{} A ``generic hyperbolic knot'' will provide 
the unique knot complement in its commensurability class.  More precisely:

\begin{proposition}
\label{genericworks}
Let $K$ be a hyperbolic knot. If $K$ admits 
no symmetries, no hidden symmetries, 
and no lens space surgeries, then $S^3\setminus K$
is the only knot complement in
its commensurability class.\end{proposition}

\noindent{\bf Proof:}~Since $K$ admits no symmetries and no hidden symmetries,
$S^3\setminus K$ is non-arithmetic and will be the minimal element in its
commensurability class. Hence any other knot complement commensurable
with $S^3\setminus K$ covers $S^3\setminus K$.  This covering is cyclic
by \cite{GAW}, and Theorem \ref{cyclic} provides a lens space 
surgery which is a contradiction.\qed\\[\baselineskip]
It is conjectured (see \cite{Go}) that if a knot $K$ admits a lens
space surgery then $K$ is tunnel number one, so in particular will
admit an order 2 symmetry which is a strong involution.  
Thus conjecturally any knot
without symmetries or hidden symmetries is the only knot complement in its commensurabilty class.\\[\baselineskip]
\noindent{\bf Example:}~The knot $9_{32}$ provides an example of a hyperbolic
knot with no symmetries, no hidden symmetries and no lens space
surgeries.  Indeed, Riley \cite{Ri3} shows that this knot complement
has no symmetries and that its trace field has degree 29. Hence there are
no hidden symmetries (see \S 2.2).  The computation of the
Alexander polynomial shows there are no
lens space surgeries by \cite{OS}.

\subsection{}
Based on the above discussion, we have the following conjecture.

\begin{conjecture} 
\label{3conj}
 Let $K$ be a hyperbolic knot. Then\\[\baselineskip]
\noindent{\bf (i)}~There are at most three knot complements in the commensurability class of $S^3\setminus K$.\\[\baselineskip]
\noindent{\bf (ii)}~If $K$ does not admit symmetries or hidden symmetries then
there is only one knot complement in the commensurability class of $S^3\setminus K$. \end{conjecture}

We summarize what is known to us. 

\begin{theorem}
\label{allwecando}
Let $K$ be a hyperbolic knot in $S^3$, then Conjecture \ref{3conj}
holds for $K$ if:

\smallskip

\noindent {(i)}~$K$ is 2-bridge, or

\smallskip

\noindent{(ii)}~$K$ admits no symmetries, no hidden symmetries
and has no lens space surgery, or

\smallskip
\noindent{(iii)}~$K$ admits a free symmetry but no other symmetries and no hidden symmetries, or
\smallskip
\noindent{(iv)}~$K$ admits a strong involution but no
other symmetries and no hidden symmetries.\end{theorem}

\noindent{\bf Proof:}~By Theorem \ref{main} and Proposition
\ref{genericworks}, it remains to prove (iii) and (iv).  For (iii), since
the action is free, and there are no hidden symmetries, the
minimal orbifold in the commensurability class of $S^3\setminus K$ is a
manifold, and so this case is handled directly by \cite[Theorem
4]{Re1}.  

For case (iv), since $K$ has no hidden symmetries the minimal orbifold
$Q_K$ (in the previous notation) is ${\bf H}^3/N(K)$. By assumption,
$[N(K):\Gamma_K]=2$.  Suppose that $S^3\setminus K'={\bf
  H}^3/\Gamma_{K'}$ is commensurable with $S^3\setminus K$, where
$\Gamma_{K'}<N(K)$.  We can assume that $S^3\setminus K'$ does not
cover $S^3\setminus K$. For if so, by \cite{WW} $S^3\setminus K'$
cannot cover any other knot complement, and the result follows from
\cite{CGLS}.  Since the only symmetry of $K$ is a strong involution
(which we shall denote by $\tau$), $N(K)$ is generated by elements of
order $2$. Hence the abelianization of $N(K)$ is generated by elements
of order $2$, and it is easy to see that all elements in the
abelianization have order $2$. In particular any non-trivial cyclic
quotient of $N(K)$ has order $2$.

Since $Q_K$ does not have a rigid cusp, Proposition \ref{nr} shows
that $S^3\setminus K'$ is also a regular cover of $Q_K$.  Let $G$ be
the covering group of $S^3\setminus K' \rightarrow Q_K$.  We claim
that $G$ is cyclic. Consider the cover $M$ of $S^3\setminus K$ and
$S^3\setminus K'$ corresponding to $\Gamma_K\cap \Gamma_{K'}$.  Since
$\Gamma_{K'}$ is assumed not to be a subgroupof $\Gamma_K$, it follows
that $M\rightarrow S^3\setminus K'$ is a 2-fold cover.  Hence, $M$ has
one cusp.  Now $M\rightarrow S^3\setminus K$ is also a regular cover,
and the covering group is necessarily $G$ (by index).  The covering
group is abelian since it is determined by the action on the cusp, and
hence cyclic since it is the quotient of a knot group. As remarked
above, any cyclic quotient of $N(K)$ has order $2$. Since the cusp of
$N(K)$ is $S^2(2,2,2,2)$, there is one two-fold cover of ${\bf H}^3/N(K)$
with a torus cusp (as in the proof of Lemma \ref{oneZ4} (ii)). Hence
$K=K'$ in this case.\qed\\[\baselineskip]
We have the following Corollary of Theorem \ref{allwecando}(iii).

\begin{corollary}
\label{applywz}
If $K$ admits no hidden symmetries and has
a lens space surgery, then Conjecture \ref{3conj} holds for $K$.
\end{corollary}

\noindent{\bf Proof:}~Since $K$ admits a lens space surgery,
\cite{WZ} shows that the only possibility for a non-trivial
symmetry of $K$ is a strong involution. Thus either we can 
apply (iv) of Theorem \ref{allwecando}, or $S^3\setminus K$ is the 
minimal element in the commensurability class. In which case,
any knot complement in the commensurability class covers $S^3\setminus K$, 
and we are done by \cite{GAW}, Theorem \ref{cyclic} and
\cite{CGLS}.\qed\\[\baselineskip]
 The trace field of the knot complement of the $(-2,3,7)$-pretzel knot has degree 3. Therefore by Proposition \ref{nr} it admits no hidden symmetries and  there are exactly 3 knot complements in this commensurability class. It seems likely that there are no examples of knots that have hidden symmetries
and admit a lens space surgery.

\bigskip

\noindent Department of Mathematics,\\
University of Texas,\\
Austin, TX 78712\\

\noindent Department of Mathematics\\
Tufts University\\
Medford, MA 02155\\

\noindent and\\

\noindent 
D\'epartement de Math\'ematiques\\
Universit\'e du Qu\'ebec \`a Montr\'eal\\
Montr\'eal, Qu\'ebec H3C 3J7 \\

 \end{document}